\documentclass[a4paper,twoside,12pt]{article}
\usepackage{amssymb} 
\begin{document}
\title{Joint-tree model and the maximum genus of  graphs \footnotemark[2]
\author{Guanghua Dong$^{1,2}$, Ning Wang$^{3}$, Yuanqiu Huang$^{1}$ and Yanpei Liu$^{4}$\\
{\small\em 1.Department of Mathematics, Normal
University of Hunan, Changsha, 410081, China}\\
\hspace{-1mm}{\small\em 2.Department of Mathematics, Tianjin
Polytechnic
University, Tianjin, 300387, China}\\
\hspace{-5mm}{\small\em 3.Department of Information Science and
Technology,
Tianjin University of Finance }\\
\hspace{-74mm} {\small\em and Economics, Tianjin, 300222, China}\\
\hspace{-8mm}{\small\em 4.Department of Mathematics, Beijing
Jiaotong University, Beijing 100044, China}\\ }}
\footnotetext[2]{\footnotesize \em   This work was partially
Supported  by the China Postdoctoral Science Foundation funded
project (Grant No: 20110491248), the New Century Excellent Talents
in University (Grant No: NCET-07-0276), and the National Natural
Science Foundation of China (Grant No: 11171114).}
\footnotetext[1]{\footnotesize \em E-mail: gh.dong@163.com(G. Dong);
ninglw@163.com(N. Wang);  hyqq@hunnu.edu.cn(Y. Huang);
ypliu@bjtu.edu.cn(Y. Liu). }

\date{}
\maketitle

\vspace{-.8cm}
\begin{abstract}

The vertex $v$ of a graph $G$ is called a 1-$critical$-$vertex$ for
the maximum genus of the graph, or for simplicity called
1-$critical$-$vertex$, if $G-v$ is a connected graph and
$\gamma_{M}(G-v)=\gamma_{M}(G)-1$. In this paper, through the
$joint$-$tree$ model, we obtained some types of
1-$critical$-$vertex$, and get the upper embeddability of the
$Spiral$ $S_{m}^{n}$.

\bigskip
\noindent{\bf Key Words:} joint-tree; maximum genus; graph embedding \\
{\bf MSC(2000):} \ 05C10
\end{abstract}


\bigskip
\noindent {\bf 1. Introduction}

In 1971,  Nordhaus, Stewart and White \cite{nor} introduced the idea
of the  maximum genus of  graphs. Since then many researchers have
paid  attention to this object and obtained many interesting
results, such as the results in [2-8] \cite{ren}
\cite{sko}\cite{xuo} etc. In this paper, by means of the joint-tree
model, which is originated from the early works of Liu (\cite{liu1})
and is formally established  in \cite{liu3} and \cite{liu5}, we
offer a method which is different from others  to find the maximum
genus of some types of graphs.

Surfaces considered here are compact 2-dimensional manifolds without
boundary. An orientable surface $S$ can be regarded as a polygon
with even number of directed edges such that both $a$ and $a^{-1}$
occurs once on $S$ for each $a\in S$, where the power
\textquotedblleft $-1$\textquotedblright  means that the direction
of $a^{-1}$ is opposite to that of $a$ on the polygon. For
convenience, a polygon is represented by a linear sequence of
lowercase letters. An elementary result in algebraic topology states
that each orientable surface is equivalent to one of the following
standard forms of surfaces:
\[ O_{p}=\left\{
 \begin{array}{ll}
 a_{0}a_{0}^{-1},  &\mbox{$p = 0$,}\\
 \prod\limits_{i=1}^{p}a_{i}b_{i}a_{i}^{-1}b_{i}^{-1},  &\mbox{$p \geq 1$ .}
 \end{array}
 \right.
\]
which are the sphere ($p=0$), torus ($p=1$), and the orientable
surfaces of genus $p\ (p\geq2)$. The genus of a surface $S$ is
denoted by $g(S)$. Let  $A$, $B$, $C$, $D$, and $E$ be possibly
empty linear sequence of letters. Suppose $A=a_1a_2\dots a_{r},
r\geq1$, then $A^{-1}=a_{r}^{-1}\dots a_2^{-1}a_1^{-1}$ is called
the $inverse$ of $A$. If  $\{a, b, a^{-1}, b^{-1}\}$ appear in a
sequence with the form as $AaBbCa^{-1}Db^{-1}E$, then they are said
to be an $interlaced$ $set$; otherwise, a $parallel$ $set$. Let
$\widetilde{S}$ be the set of all surfaces. For a surface $S\in
\widetilde{S}$, we obtain its genus $g(S)$ by using the following
transforms to determine its equivalence to one of the standard
forms.

\medskip

{\bf Transform 1} \ \ $Aaa^{-1} \sim A$, where $A\in \widetilde{S}$
and $a\notin A$.

{\bf Transform 2} \ \  $AabBb^{-1}a^{-1} \sim AcBc^{-1}$.

{\bf Transform 3} \  \ $(Aa)(a^{-1}B) \sim (AB)$.

{\bf Transform 4} \ \  $AaBbCa^{-1}Db^{-1}E \sim
ADCBEaba^{-1}b^{-1}$.

\medskip

\noindent   In the above transforms, the parentheses stand for
cyclic order. For convenience, the parentheses are always omitted
when unnecessary to distinguish cyclic or linear order. For more
details concerning surfaces, the reader is referred to \cite{liu3},
\cite{liu5} and \cite{rin}.

Let $T$ be a spanning tree of a graph $G=(V,E)$, then
$E=E_{T}+E^{*}_{T}$, where $E_{T}$ consists of all the tree edges,
and $E^{*}_{T}=\{ e_1, e_2, \dots e_{\beta}\}$ consists of all the
co-tree edges, where $\beta=\beta(G)$ is the cycle rank of $G$.
Split each co-tree edge $e_{i}=(\mu_{e_{i}}, \nu_{e_{i}}) \in
E^{*}_{T}$ into two semi-edges $(\mu_{e_{i}}, \omega_{e_{i}})$,
$(\nu_{e_{i}}, \omega'_{e_{i}})$, denoted by $e_{i}^{+1}$ (or simply
by $e_{i}$ if no confusion) and $e_{i}^{-1}$ respectively.  Let
$\widetilde{T}=(V+V_1, E+E_1)$, where $V_1=\{\omega_{e_{i}}, \
\omega'_{e_{i}} \ | \ 1\leqslant i \leqslant \beta \}$,
$E_1=\{(\mu_{e_{i}}, \omega_{e_{i}}), \ (\nu_{e_{i}},
\omega'_{e_{i}}) \  | \ 1\leqslant i \leqslant \beta\}$. Obviously,
$\widetilde{T}$ is a tree. A $rotation$ $at$ $a$ $vertex$ $v$, which
is denoted by $\sigma_{v}$, is a cyclic permutation of edges
incident on $v$. A rotation system $\sigma=\sigma_{G}$ for a graph
$G$  is a set $\{\sigma_{v} | \forall v \in V(G)\}$. The tree
$\widetilde{T}$ with a rotation system of $G$ is called a
$joint$-$tree$ of $G$, and is denoted by $\widetilde{T}_{\sigma}$.
Because it ia a tree, it can be embedded in the plane. By reading
the lettered semi-edges of $\widetilde{T}_{\sigma}$ in a fixed
direction (clockwise or anticlockwise), we can get an algebraic
representation of the surface which is represented by a
$2\beta-$polygon. Such a surface, which is denoted by $S_{\sigma}$,
is called an associated surface of $\widetilde{T}_{\sigma}$. A
joint-tree $\widetilde{T}_{\sigma}$ of $G$ and its associated
surface is illustrated by Fig.1, where the rotation at each vertex
of $G$ complies with the clockwise rotation. From \cite{liu3}, there
is 1-1 correspondence between associated surfaces (or joint-trees)
and embeddings of a graph.


\bigskip
\setlength{\unitlength}{1mm}
\begin{center}
\begin{picture}(100,40)

\put(15,3){\begin{picture}(10,10)

\put(29,-5){{\bf Fig. 1.}}
\end{picture}
}


\put(-15,4){\begin{picture}(10,10)

\put(0,10){\circle*{1.8}}

\put(26,10){\circle*{1.8}}

\put(13,28){\circle*{1.8}}

\put(13,17){\circle*{1.8}}

\put(1,10){\line(1,0){24}}

\thicklines

\put(0.8,10.5){\line(2,1){11.5}}

\put(25,10.5){\line(-2,1){11.2}}

\put(13,27){\line(0,-1){9}}

\thinlines

\put(0.7,11){\line(3,4){12}}

\put(25.7,11){\line(-3,4){12}}

\put(12,0){ $G$}

\footnotesize

{ \put(0.5,20){ $e_1$}

\put(20,20){ $e_2$}

\put(9,6){ $e_3$}  }

\end{picture}
}

\put(35,4){\begin{picture}(10,10)

\put(7,14){\circle*{1.8}}

\put(19,14){\circle*{1.8}}

\put(13,22.5){\circle*{1.8}}

\put(13,17){\circle*{1.8}}

\put(-1,21){\circle*{1.8}}  

\put(-1,7){\circle*{1.8}}   

\put(6,30){\circle*{1.8}}  

\put(20,30){\circle*{1.8}}  

\put(27,21){\circle*{1.8}}  

\put(27,7){\circle*{1.8}}  

\thicklines

\put(13,17){\line(-2,-1){6}}

\put(13,17){\line(2,-1){6}}

\put(13,17){\line(0,1){6}}

\thinlines

\put(6.5,14){\line(-1,1){7}}

\put(6.5,14){\line(-1,-1){7}}

\put(19.5,14){\line(1,1){7}}

\put(19.5,14){\line(1,-1){7}}

\put(13,23){\line(-1,1){7}}

\put(13,23){\line(1,1){7}}

\put(11,2){ $\widetilde{T}_{\sigma}$}

\footnotesize

\put(-0.9,16){ $e_1$}

\put(4,22){ $e_1^{-1}$}

\put(12,27){ $e_2$}

\put(17,18){ $e_2^{-1}$}

\put(22,11){ $e_3$}

\put(2,8){ $e_3^{-1}$}

\put(-5,23){ $\omega_{e_1}$}

\put(-8,7){ $\omega^{'}_{e_3}$}

\put(4,32){ $\omega^{'}_{e_1}$}

\put(18,32){ $\omega_{e_2}$}

\put(27,21){ $\omega^{'}_{e_2}$}

\put(27,7){ $\omega_{e_3}$}

\end{picture}
}

\put(85,4){\begin{picture}(10,10)

\put(3,28){\line(1,0){15}}

\put(3,28){\line(-3,-5){5.5}}

\put(3,10){\line(1,0){15}}

\put(3,10){\line(-3,5){5.5}}

\put(18,10){\line(3,5){5.5}}

\put(18,28){\line(3,-5){5.5}}

\put(12,0){ $S_{\sigma}$}

\footnotesize

 \put(-5,24){ $e_1$}

\put(7,30){ $e_1^{-1}$}

\put(21,23){ $e_2$}

\put(21,13){ $e_2^{-1}$}

\put(8,5){ $e_3$}

\put(-7,11){ $e_3^{-1}$}

\put(8,18){ $\curvearrowright$ }

\end{picture}
}


\end{picture}
\end{center}

To $merge$ a vertex of degree two is that replace its two incident
edges with a single edge joining the other two incident vertices.
$Vertex$-$splitting$ is such an operation as follows. Let $v$ be a
vertex of graph $G$. We replace $v$ by two new vertices $v_1$ and
$v_2$. Each edge of $G$ joining $v$ to another vertex $u$ is
replaced by an edge joining $u$ and $v_1$, or by an edge joining $u$
and $v_2$.  A graph is called a $cactus$ if all circuits are
independent, $i.e.$, pairwise vertex-disjoint. The $maximum$ $genus$
$\gamma_M(G)$ of a connected graph \emph{G} is the maximum integer
\emph{k} such that there exists an embedding of $G$ into the
orientable surface of genus $k$. Since any embedding must have at
least one face, the Euler characteristic for one face leads to an
upper bound on the maximum genus
\begin{displaymath}
\gamma_M(G)\leq \lfloor \frac{|E(G)| - |V(G)| + 1}{2} \rfloor.
\end{displaymath}
A graph $G$ is said to be $upper$ $embeddable$ if $\gamma_M(
\emph{G} )$ = $\lfloor\frac{\beta(G)}{2}\rfloor$, where $\beta(
\emph{G} )=| E(G)|$ $-$ $|V(G)|$ + 1 denotes the \emph{Betti number}
of  \emph{G}. Obviously, the maximum genus of a cactus is zero.
 The vertex $v$ of a graph $G$ is called a 1-$critical$-$vertex$ for
the maximum genus of the graph, or for simplicity called
1-$critical$-$vertex$, if $G-v$ is a connected graph and
$\gamma_{M}(G-v)=\gamma_{M}(G)-1$. Graphs considered here are all
connected, undirected, and with minimum degree at least three. In
addition, the surfaces are all orientable. Notations and
terminologies not defined here can be seen in \cite{bon},
\cite{liu4}, \cite{liu3},  and \cite{liu5}.

\bigskip

{\bf Lemma 1.0}   \ \ If there is a joint-tree
$\widetilde{T}_{\sigma}$ of $G$ such that the genus of its
associated surface equals $\lfloor\frac{\beta(G)}{2}\rfloor$ then
$G$ is upper embeddable.

\bigskip

{\bf Proof } \ \ According to the definition of joint-tree,
associated surface, and upper embeddable graph, Lemma 1.0 can be
easily obtained.  $\hspace*{\fill} \Box$

\bigskip

{\bf Lemma 1.1}   \ \ Let $AB$ be a surface. If $x\notin A\cup B$,
then $g(AxBx^{-1})=g(AB)$ or $g(AxBx^{-1})=g(AB)+1$.

\bigskip

{\bf Proof } \ \ First discuss the topological standard form of the
surface $AB$. {\bf (I)} According to the left to right direction,
let $\{x_1, y_1, x_1^{-1}, y_1^{-1}\}$ be the first interlaced set
appeared in $A$. Performing Transform 4 on $\{x_1, y_1, x_1^{-1},
y_1^{-1}\}$ we will get $A^{'}Bx_1y_1x_1^{-1}y_1^{-1}$ ($\sim$
$AB$). Then perform Transform 4 on the first interlaced set in
$A^{'}$. And so on. Eventually we will get
$\widetilde{A}B\prod\limits_{i=1}^{r}x_{i}y_{i}x_{i}^{-1}y_{i}^{-1}$
($\sim$ $AB$), where there is no interlaced set in $\widetilde{A}$.
{\bf (II)} For the surface
$\widetilde{A}B\prod\limits_{i=1}^{r}x_{i}y_{i}x_{i}^{-1}y_{i}^{-1}$,
from the left of $B$, successively perform Transform 4 on $B$
similar to that on $A$ in (I). Eventually we will get
$\widetilde{A}\widetilde{B}\prod\limits_{i=1}^{r}x_{i}y_{i}x_{i}^{-1}y_{i}^{-1}
\prod\limits_{j=1}^{s}a_{j}b_{j}a_{j}^{-1}b_{j}^{-1}$ ($\sim$ $AB$),
where there is no interlaced set in $\widetilde{B}$. {\bf (III)} For
the surface
$\widetilde{A}\widetilde{B}\prod\limits_{i=1}^{r}x_{i}y_{i}x_{i}^{-1}y_{i}^{-1}
\prod\limits_{j=1}^{s}a_{j}b_{j}a_{j}^{-1}b_{j}^{-1}$, from the left
of $\widetilde{A}\widetilde{B}$, successively perform Transform 4 on
$\widetilde{A}\widetilde{B}$ similar to that on $A$ in (I). At last,
we will get $\prod\limits_{i=1}^{p}a_{i}b_{i}a_{i}^{-1}b_{i}^{-1}$,
which is the topologically standard form of the surface $AB$.

As for the surface $AxBx^{-1}$, perform Transform 4 on $A$ and $B$
similar to that on $A$ in (I) and  $B$ in (II) respectively.
Eventually
$\widetilde{A}x\widetilde{B}x^{-1}\prod\limits_{i=1}^{r}x_{i}y_{i}x_{i}^{-1}y_{i}^{-1}
\prod\limits_{j=1}^{s}a_{j}b_{j}a_{j}^{-1}b_{j}^{-1}$ ($\sim$
$AxBx^{-1}$) will be obtained. Then perform  the same Transform 4 on
$\widetilde{A}x\widetilde{B}x^{-1}$  as that on
$\widetilde{A}\widetilde{B}$ in (III), and at last, one more
Transform 4 than that in (III) may be needed  because of $x$ and
$x^{-1}$ in $\widetilde{A}x\widetilde{B}x^{-1}$. Eventually
$\prod\limits_{i=1}^{p}a_{i}b_{i}a_{i}^{-1}b_{i}^{-1}$ or
$\prod\limits_{i=1}^{p+1}a_{i}b_{i}a_{i}^{-1}b_{i}^{-1}$, which is
the topologically standard form of the surface $AxBx^{-1}$, will be
obtained.

From the above, Lemma 1.1 is obtained. $\hspace*{\fill} \Box$

\bigskip

 {\bf Lemma 1.2}   \ \  Among all orientable surfaces
represented by the linear sequence consisting of $a_{i}$ and
$a_{i}^{-1}$ ($i=1, \dots, n$), the surface $a_1a_2\dots
a_{n}a_1^{-1}a_2^{-1}\dots a_{n}^{-1}$ is   one whose genus is
maximum.

\bigskip

{\bf Proof } \ \ According to Transform 4, Lemma 1.2 can be easily
obtained. $\hspace*{\fill} \Box$

\bigskip

{\bf Lemma 1.3}   \ \ Let $G$ be a graph with minimum degree at
least three, and $\bar{G}$ be the graph obtained from $G$ by a
sequence of vertex-splitting,  then $\gamma_{M}(\bar{G}) \leq
\gamma_{M}(G)$. Furthermore, if $\bar{G}$ is upper embeddable then
$G$ is upper embeddable as well.

\bigskip

{\bf Proof } \ \ Let $v$ be a vertex of degree $n(\geq4)$ in $G$,
and $G^{'}$ be the graph obtained from $G$ by splitting the vertex
$v$ into two vertices such that both their degrees are at least
three. First of all, we prove that the maximum genus will not
increase after one vertex-splitting operation, $i.e.$,
$\gamma_{M}(G^{'}) \leq \gamma_{M}(G)$.

 Let $e_1$, $e_2$, $\dots$ $e_{n}$ be the $n$ edges
incident to $v$, and $v$ be split into $v_1$ and  $v_2$. Without
loss of generality, let $e_{i_1}$, $e_{i_2}$, $\dots$ $e_{i_{r}}$ be
incident to $v_1$, and $e_{i_{r+1}}$,  $\dots$ $e_{i_{n}}$ be
incident to $v_2$, where $2\leq i_{r}\leq n-2$. Select such a
spanning tree $T$ of $G$ that $e_{i_1}$ is a tree edge, and
$e_{i_2}$, $\dots$ $e_{i_{n}}$ are all co-tree edges. As for graph
$G^{'}$, select $T^{*}$ be a spanning tree such that both  $e_{i_1}$
and $(v_1, v_2)$ are tree edges, and the other edges of $T^{*}$ are
the same as the edges in $T$. Obviously,  $e_{i_2}$, $\dots$
$e_{i_{n}}$ are  co-tree edges of $T^{*}$. Let
$\mathcal{T}$=$\{\hat{T}_{\sigma}|\hat{T}_{\sigma}=\overline{(T-v)}_{\sigma},$
where $\overline{(T-v)}_{\sigma}$ is a joint-tree of $G-v$\},
$\mathcal{T}^{*}$=$\{\hat{T}^{*}_{\sigma}|\hat{T}^{*}_{\sigma}=\overline{(T^{*}-\{v_1,
v_2\})}_{\sigma},$ where $\overline{(T^{*}-\{v_1, v_2\})}_{\sigma}$
is a joint-tree of $G^{'}-\{v_1, v_2\}$\}. It is obvious that
$\mathcal{T}=\mathcal{T}^{*}$.  Let $\mathcal {S}$ be the set of all
the associated surfaces of the joint-trees of $G$, and $\mathcal
{S}^{*}$ be the set of all the associated surfaces of the joint
trees of $G^{'}$. Obviously, $\mathcal {S}^{*} \subseteq \mathcal
{S}$.  Furthermore, $|\mathcal {S}^{*}|=r!\times
(n-r)!\times|\mathcal{T}^{*}|$ $<$ $|\mathcal
{S}|=(n-1)!\times|\mathcal{T}|$. So $\mathcal {S}^{*} \subset
\mathcal {S}$, and  we have $\gamma_{M}(G^{'}) \leq \gamma_{M}(G)$.

Reiterating this procedure, we can get that $\gamma_{M}(\bar{G})
\leq \gamma_{M}(G)$. Furthermore, because $\beta(G)=\beta(\bar{G})$,
it can be obtained that if $\bar{G}$ is upper embeddable then
$\lfloor\frac{\beta(G)}{2}\rfloor$ =
$\lfloor\frac{\beta(\bar{G})}{2}\rfloor$ = $\gamma_M( \bar{G})$
$\leq \gamma_M(G)$ $\leq \lfloor\frac{\beta(G)}{2}\rfloor$. So,
$\gamma_M(G) = \lfloor\frac{\beta(G)}{2}\rfloor$, and $G$ is upper
embeddable. $\hspace*{\fill} \Box$

\bigskip

\bigskip

\noindent {\bf 2. Results related to 1-critical-vertex }

\bigskip

The $neckband$ $\mathcal {N}_{2n}$ is such a graph that $\mathcal
{N}_{2n}=C_{2n}+R$, where $C_{2n}$ is a 2n-cycle, and $R=\{ a_{i} |
a_{i}=(v_{2i-1}, v_{2i+2}).\ (i=1,2,\dots, n, \ 2i+2\equiv r (mod\
2n),\ 1\leq r < 2n)\}$. The $m\ddot{o}bius$ $ladder$ $\mathcal
{M}_{2n}$ is such  a cubic circulant graph with 2n vertices, formed
from a 2n-cycle by adding edges (called "rungs") connecting opposite
pairs of vertices in the cycle. For example, Fig. 2.1 and Fig. 2.5
is a graph of $\mathcal {N}_8$ and $\mathcal {M}_{2n}$ respectively.
A vertex like the solid vertex in Fig. 2.2, Fig. 2.3, Fig. 2.4, Fig.
2.5, and Fig. 2.6 is called an $\alpha$-$vertex$, $\beta$-$vertex$,
$\gamma$-$vertex$, $\delta$-$vertex$, and $\eta$-$vertex$
respectively, where Fig. 2.6 is a neckband.

\begin{footnotesize}
\setlength{\unitlength}{1mm}
\begin{center}
\begin{picture}(100,30)

\put(-13,-13){\begin{picture}(10,10)

\put(10,30){\circle{15}}

\put(6.8,36.4){\circle*{1.5}}

\put(4,39){$v_3$}

\put(13.4,36.4){\circle*{1.5}}

\put(13,38.5){$v_4$}

\put(3.2,32){\circle*{1.5}}

\put(-1.8,32){$v_2$}

\put(16.7,32){\circle*{1.5}}

\put(18.5,32){$v_5$}

\put(3.7,27){\circle*{1.5}}

\put(-1,25.5){$v_1$}

\put(16.3,27){\circle*{1.5}}

\put(18,25.5){$v_6$}

\put(13.4,24){\circle*{1.5}}

\put(13,20){$v_7$}

\put(6.8,24){\circle*{1.5}}

\put(5,20){$v_8$}

\qbezier(3.7,27)(6,36)(13.4,36.4)

\qbezier(6.8,36.4)(14,36)(16.3,27)

\qbezier(16.7,32)(13.5,24)(6.8,24)

\qbezier(3.2,32)(6,24)(13.4,24)

\put(3,14){{\bf Fig. 2.1.}}
\end{picture}}

\put(40,-13){\begin{picture}(10,10)

\thicklines

\put(0,30){\circle{6}}

\thinlines

\put(2.8,30){\circle*{1.5}}

\put(11.2,30){\circle{1.5}}

\put(2.5,30){\line(1,0){8}}

\put(11.5,30.6){\line(1,1){6}}

\put(11.5,29.4){\line(1,-1){6}}

\put(16,32){\circle*{0.8}}

\put(16,30){\circle*{0.8}}

\put(16,28){\circle*{0.8}}

\put(-4.5,32){$a$}

\put(2,14){{\bf Fig. 2.2}}
\end{picture}}


\put(85,-13){\begin{picture}(10,10)

\put(8,30){\circle*{1.5}}

\put(16,30){\circle*{1.5}}

\thicklines

\qbezier(8,30)(12,34)(16,30)

\qbezier(8,30)(12,26)(16,30)

\thinlines

\put(1,30){\circle{1.5}}

\put(23,30){\circle{1.5}}

\put(8,30){\line(-1,0){6.5}}

\put(16,30){\line(1,0){6.3}}

\put(0.5,30.5){\line(-1,1){5}}

\put(0.5,29.5){\line(-1,-1){5}}

\put(23.5,30.5){\line(1,1){5}}

\put(23.5,29.5){\line(1,-1){5}}

\put(-3.5,31.8){\circle*{0.8}}

\put(-3.5,29.8){\circle*{0.8}}

\put(-3.5,27.8){\circle*{0.8}}

\put(27.5,31.8){\circle*{0.8}}

\put(27.5,29.8){\circle*{0.8}}

\put(27.5,27.8){\circle*{0.8}}

\put(7,14){{\bf Fig. 2.3}}

\put(12,33){$a$}

\put(12,24){$b$}
\end{picture}}

\end{picture}
\end{center}

\setlength{\unitlength}{1mm}
\begin{center}
\begin{picture}(100,34)

\put(-13,-10){\begin{picture}(10,10)

\put(7,30){\circle{1.5}}

\put(17,30){\circle{1.5}}

\put(12,34){\circle*{1.5}}

\put(12,26){\circle*{1.5}}

\put(12,34){\line(-5,-4){4.4}}

\thicklines

\put(12,34){\line(5,-4){4.4}}

\put(12,34){\line(0,-1){8}}

\thinlines

\thinlines

 \put(12,26){\line(-5,4){4.4}}

\put(12,26){\line(5,4){4.4}}

\put(1,30){\circle{1.5}}

\put(23,30){\circle{1.5}}

\put(6.5,30){\line(-1,0){5}}

\put(17.5,30){\line(1,0){5}}

\put(0.5,30.5){\line(-1,1){5}}

\put(0.5,29.5){\line(-1,-1){5}}

\put(23.5,30.5){\line(1,1){5}}

\put(23.5,29.5){\line(1,-1){5}}

\put(-3.5,31.8){\circle*{0.8}}

\put(-3.5,29.8){\circle*{0.8}}

\put(-3.5,27.8){\circle*{0.8}}

\put(27.5,31.8){\circle*{0.8}}

\put(27.5,29.8){\circle*{0.8}}

\put(27.5,27.8){\circle*{0.8}}

\put(10.6,35.8){$v_1$}

\put(11,22.3){$v_2$}

\put(14.5,32.5){$a$}

\put(10,28.5){$b$}

\put(6,13){{\bf Fig. 2.4}}
\end{picture}}


\put(29,-10){\begin{picture}(10,10)

\thinlines

 \put(20,30){\circle{15}}

\put(15.5,25){\circle*{1.5}}

\put(15.5,35){\circle*{1.5}}

\put(20,37){\circle*{1.5}}

\put(20,23){\circle*{1.5}}

\put(25,25){\circle*{1.5}}

\put(25,35){\circle*{1.5}}

\put(20,37){\line(0,-1){14}}

\thicklines

\put(15.2,25){\line(1,1){9.5}}

\put(15,35.2){\line(1,-1){9.5}}

\qbezier(15.5,35)(16.8,37)(20,37)

\qbezier(20,37)(22.5,37)(25,35)

\thinlines

\put(24.5,32){\circle*{0.8}}

\put(24.5,30){\circle*{0.8}}

\put(24.5,28){\circle*{0.8}}

\put(15.5,32){\circle*{0.8}}

\put(15.5,30){\circle*{0.8}}

\put(15.5,28){\circle*{0.8}}

\put(19,40){$v_1$}

\put(27,35){$v_2$}

\put(27,23){$v_{n}$}

\put(18,19){$v_{n+1}$}

\put(8,22){$v_{n+2}$}

\put(8,35){$v_{2n}$}

\put(14.5,37){$m$}

\put(22.5,37){$n$}

\put(16,13){{\bf Fig. 2.5}}
\end{picture}}


\put(78,-10){\begin{picture}(10,10)

\thinlines

\put(20,30){\circle{15}}

\put(17.5,36.5){\circle*{1.5}}

\put(15,38.5){$v_2$}

\put(22.5,36.5){\circle*{1.5}}

\put(20,39){$v_1$}

\put(14.5,34.2){\circle*{1.5}}

\put(10,35){$v_3$}

\put(25.5,34.5){\circle*{1.5}}

\put(27.5,35.5){$v_{2n}$}

\put(13,30.5){\circle*{1.5}}

\put(8.2,30){$v_4$}

\put(27,30.5){\circle*{1.5}}

\put(28.5,29){$v_{2n-1}$}

\put(25.7,26){\circle*{1.5}}

\put(27,24){$v_{2n-2}$}

\put(22.2,23.3){\circle*{1.5}}

\put(21,20){$v_{2n-3}$}

\thicklines

\qbezier(22.5,36.5)(24,36)(25.5,34.5)

\qbezier(27,30.5)(26.3,34)(25.5,34.5)

\qbezier(22.5,36.5)(19,30)(13,30.5)

\qbezier(17.5,36.5)(21,30)(27,30.5)

\multiput(14.5,34.2)(0.2,-1.6){5}{\circle*{0.6}}

\multiput(25.7,26)(-1.5,-0.3){5}{\circle*{0.6}}

\thinlines

\qbezier(25.5,34.5)(21,29.5)(22.2,23.3)

\put(24.5,36){$s$}

\put(27.3,32){$r$}

\put(15,13){{\bf Fig. 2.6}}

\end{picture}}

\end{picture}
\end{center}
\end{footnotesize}

\bigskip

{\bf Theorem 2.1}   \ \  If $v$ is an $\alpha$-$vertex$ of a graph
$G$, then $\gamma_{M}(G-v)$ = $\gamma_{M}(G)$. If $v$ is a
$\beta$-$vertex$, or a $\gamma$-$vertex$, or a $\delta$-$vertex$, or
an $\eta$-$vertex$  of a graph $G$, and $G-v$ is a connected graph,
then $\gamma_{M}(G-v)$ = $\gamma_{M}(G)-1$, $i.e.$,
$\beta$-$vertex$, $\gamma$-$vertex$, $\delta$-$vertex$ and
$\eta$-$vertex$ are 1-$critical$-$vertex$.

\bigskip

{\bf Proof } \ \ If $v$ is an $\alpha$-$vertex$ of the graph $G$,
then it is easy to get that $\gamma_{M}(G-v)$ = $\gamma_{M}(G)$. In
the following, we will discuss the other cases.

$\textbf{Case 1:}$ \ $v$ is an $\beta$-$vertex$ of $G$.

According to Fig. 2.3, select such a spanning tree $T$ of $G$ such
that both $a$ and $b$ are co-tree edges. It is obvious that the
associated surface for each joint-tree of $G$ must be one of the
following four forms: (i) $AabBa^{-1}b^{-1}$ $\sim
ABaba^{-1}b^{-1}$, (ii) $AabBb^{-1}a^{-1}$ $\sim AcBc^{-1}$, (iii)
$AbaBa^{-1}b^{-1}$ $\sim AcBc^{-1}$, (iv) $AbaBb^{-1}a^{-1}$ $\sim
ABbab^{-1}a^{-1}$. On the other hand, for each joint-tree
$\widetilde{T}^{*}_{\sigma}$, which is a joint-tree of $G-v$, its
associated surface must be the form as $AB$, where $A$ and $B$ are
the same as that in the above four forms. According to (i)-(iv),
Lemma 1.1, and $g(ABaba^{-1}b^{-1})$=$g(AB)+1$,  we can get that
$\gamma_{M}(G-v)$ = $\gamma_{M}(G)-1$.

$\textbf{Case 2:}$ \ $v$ is an $\gamma$-$vertex$ of $G$.

As illustrated by Fig. 2.4, both $v_1$ and $v_2$ are
$\gamma$-$vertex$. Without loss of generality, we only prove that
$\gamma_{M}(G-v_1)$ = $\gamma_{M}(G)-1$. Select such a spanning tree
$T$ of $G$ such that both $a$ and $b$ are co-tree edges. The
associated surface for each joint-tree of $G$ must be one of the
following 16 forms:

$$
\begin{array} {cccc}
Aabb^{-1}a^{-1}B, & Aabb^{-1}Ba^{-1}, & Aaba^{-1}Bb^{-1}, & AabBa^{-1}b^{-1},\\
Abab^{-1}a^{-1}B, & Abab^{-1}Ba^{-1}, & Abaa^{-1}Bb^{-1}, &
AbaBa^{-1}b^{-1},\\
Ab^{-1}a^{-1}Bab, & Ab^{-1}Ba^{-1}ab, & Aa^{-1}Bb^{-1}ab, & ABa^{-1}b^{-1}ab, \\
Ab^{-1}a^{-1}Bba, & Ab^{-1}Ba^{-1}ba, & Aa^{-1}Bb^{-1}ba, &
ABa^{-1}b^{-1}ba.
\end{array}
$$
Furthermore, each of these 16 types of surfaces is topologically
equivalent to one of such surfaces as $AB$, $ABaba^{-1}b^{-1}$, and
$AcBc^{-1}$. On the other hand, for each joint-tree
$\widetilde{T}^{*}_{\sigma}$, which is a joint-tree of $G-v_1$, its
associated surface must be the form of $AB$, where $A$ and $B$ are
the same as that in the above 16 forms. According to Lemma 1.1 and
$g(ABaba^{-1}b^{-1})$=$g(AB)+1$, we can get that $\gamma_{M}(G-v)$ =
$\gamma_{M}(G)-1$.

$\textbf{Case 3:}$ \ $v$ is an $\delta$-$vertex$ of $G$.

In Fig. 2.5, let $a_{i}=(v_{i}, v_{n+i}), i=1,2,\dots, n.$ Without
loss of generality, we only  prove that $\gamma_{M}(G-v_1)$ =
$\gamma_{M}(G)-1$. Select such a joint-tree $\widetilde{T}_{\sigma}$
of Fig. 2.5, which is  illustrated by Fig.3, where the edges of the
spanning tree are represented by solid line. It is obvious that the
associated surface of $\widetilde{T}_{\sigma}$ is
$mnm^{-1}n^{-1}a_2a_3\dots a_{n}a_2^{-1}a_3^{-1}\dots a_{n}^{-1}$.
On the other hand, $a_2a_3\dots a_{n}a_2^{-1}a_3^{-1}\dots
a_{n}^{-1}$ is the associated surface of one of the joint-trees of
$G-v_1$. From Lemma 1.2 and $g(mnm^{-1}n^{-1}a_2a_3\dots
a_{n}a_2^{-1}a_3^{-1}\dots a_{n}^{-1})$=$g(a_2a_3\dots
a_{n}a_2^{-1}a_3^{-1}\dots a_{n}^{-1})+1$, we can get that
$\gamma_{M}(G-v)$ = $\gamma_{M}(G)-1$.

\begin{footnotesize}
\bigskip
\setlength{\unitlength}{1mm}
\begin{center}
\begin{picture}(100,35)

\put(-13,-10){\begin{picture}(10,10)

\put(-2,30){\circle*{1.5}}

\put(3.5,30){\circle*{1.5}}

\put(6,33){\circle*{0.8}}

\put(9,33){\circle*{0.8}}

\put(12,33){\circle*{0.8}}

\put(14.5,30){\circle*{1.5}}

\put(20,30){\circle*{1.5}}

\put(25.5,30){\circle*{1.5}}

\put(27.5,33){\circle*{0.8}}

\put(30,33){\circle*{0.8}}

\put(32,33){\circle*{0.8}}

\put(34.5,30){\circle*{1.5}}

\put(40,30){\circle*{1.5}}

\put(20,25){\circle*{1.5}}

\put(-2,30){\line(1,0){42}}

\multiput(-2,30)(-1,1){5}{\circle*{0.5}}

\multiput(-2,30)(-1,-1){5}{\circle*{0.5}}

\multiput(3.5,30)(0,1){6}{\circle*{0.5}}

\multiput(14.5,30)(0,1){6}{\circle*{0.5}}

\multiput(25.5,30)(0,1){6}{\circle*{0.5}}

\multiput(34.5,30)(0,1){6}{\circle*{0.5}}

\multiput(40,30)(1,1){5}{\circle*{0.5}}

\multiput(40,30)(1,-1){5}{\circle*{0.5}}

\multiput(20,25)(1,-0.6){6}{\circle*{0.5}}

\multiput(20,25)(-1,-0.6){6}{\circle*{0.5}}

\put(20,30){\line(0,-1){5}}

\put(-3,27){$v_{2n}$}

\put(17,32){$v_{n+1}$}

\put(21,25){$v_{1}$}

\put(38,27){$v_2$}

\put(44,35){$a_2$}

\put(44,23){$n^{-1}$}

\put(33,37){$a_3$}

\put(24,37){$a_{n}$}

\put(12,37){$a_2^{-1}$}

\put(0,37){$a_{n-1}^{-1}$}

\put(-11,35){$a_{n}^{-1}$}

\put(-11,25){$m$}

\put(11.5,20.5){$n$}

\put(27,20.5){$m^{-1}$}

\put(14,10){{\bf Fig. 3.}}
\end{picture}}


\put(75,-10){\begin{picture}(10,10) \put(-2,30){\circle*{1.5}}

\put(2,30){\circle*{1.5}}

\put(7,30){\circle*{1.5}}

\put(12,30){\circle*{1.5}}

\put(17,30){\circle*{1.5}}

\put(22,30){\circle*{1.5}}

\put(24,33){\circle*{0.8}}

\put(26,33){\circle*{0.8}}

\put(28,33){\circle*{0.8}}

\put(30,30){\circle*{1.5}}

\put(35,30){\circle*{1.5}}

\put(40,30){\circle*{1.5}}

\put(22,25){\circle*{1.5}}

\put(-2,30){\line(1,0){42}}

\multiput(-2,30)(-1,1){5}{\circle*{0.5}}

\multiput(-2,30)(-1,-1){5}{\circle*{0.5}}

\multiput(2,30)(0,1){6}{\circle*{0.5}}

\multiput(7,30)(0,1){6}{\circle*{0.5}}

\multiput(12,30)(0,1){6}{\circle*{0.5}}

\multiput(17,30)(0,1){6}{\circle*{0.5}}

\multiput(22,30)(0,1){6}{\circle*{0.5}}

\multiput(35,30)(0,1){6}{\circle*{0.5}}

\multiput(40,30)(0.6,1){5}{\circle*{0.5}}

\multiput(40,30)(0.6,-1){5}{\circle*{0.5}}

\multiput(22,25)(1,-0.6){6}{\circle*{0.5}}

\multiput(22,25)(-1,-0.6){6}{\circle*{0.5}}

\put(30,30){\line(-3,-2){8}}

\put(-7,29.5){$v_{1}$}

\put(1.5,26.5){$v_{2}$}

\put(28,26){$v_{2n-3}$}

\put(20.4,21){$v_{2n}$}

\put(42,29){$v_{2n-1}$}

\put(42.7,36){$a_{n}^{-1}$}

\put(43.5,23){$r$}

\put(15.8,37){$a_3$}

\put(5.6,37){$a_2$}

\put(10.2,38){$a_1^{-1}$}

\put(21,37){$a_2^{-1}$}

\put(33,37){$a_{n-2}^{-1}$}

\put(-1,37){$a_{n}$}

\put(-11,35){$a_1$}

\put(-9,22){$s^{-1}$}

\put(10,19){$r^{-1}$}

\put(29,20.5){$s$}

\put(14,10){{\bf Fig. 4.}}
\end{picture}}


\end{picture}
\end{center}
\end{footnotesize}

$\textbf{Case 4:}$ \ $v$ is an $\eta$-$vertex$ of $G$.

As illustrated by Fig. 2.6, every vertex in Fig. 2.6 is a
$\eta$-$vertex$. Without loss of generality, we only  prove that
$\gamma_{M}(G-v_{2n})$ = $\gamma_{M}(G)-1$.

A joint-tree $\widetilde{T}_{\sigma}$ of Fig. 2.6 is depicted by
Fig.4. It can be read from Fig.4 that the associated surface of
$\widetilde{T}_{\sigma}$ is
$S=a_1a_{n}(\prod\limits_{i=1}^{n-3}a_{i+1}a_{i}^{-1})a_{n-2}^{-1}a_{n}^{-1}rsr^{-1}s^{-1}$.
Performing a sequence of Transform 4 on $S$, we have
\begin{eqnarray}
S  & = & a_1a_{n}(\prod\limits_{i=1}^{n-3}a_{i+1}a_{i}^{-1})a_{n-2}^{-1}a_{n}^{-1}rsr^{-1}s^{-1} \nonumber \\
        \mbox{(Transform\ 4)}  & \sim  &
        (\prod\limits_{i=2}^{n-3}a_{i+1}a_{i}^{-1})a_{n-2}^{-1}a_2rsr^{-1}s^{-1}a_1a_{n}a_1^{-1}a_{n}^{-1} \nonumber \\
        \mbox{(Transform\ 4)}  & \sim  &
        (\prod\limits_{i=4}^{n-3}a_{i+1}a_{i}^{-1})a_{n-2}^{-1}a_4rsr^{-1}s^{-1}a_1a_{n}a_1^{-1}a_{n}^{-1}a_3a_2a_3^{-1}a_2^{-1} \nonumber \\
        & \cdots & \ \   \cdots \nonumber \\
       \mbox{(Transform\ 4)}  & \sim  & \left\{
        \begin{array}{ll}
        rsr^{-1}s^{-1}a_1a_{n}a_1^{-1}a_{n}^{-1}(\prod\limits_{i=2}^{n-4}a_{i+1}a_{i}a_{i+1}^{-1}a_{i}^{-1}) & \  \mbox{$n \equiv 0(mod \ 2)$;} \\
        rsr^{-1}s^{-1}a_1a_{n}a_1^{-1}a_{n}^{-1}(\prod\limits_{i=2}^{n-3}a_{i+1}a_{i}a_{i+1}^{-1}a_{i}^{-1}) & \ \mbox{$n \equiv 1(mod \ 2)$.}  \\
        \end{array}
        \right.
\end{eqnarray}
It is known from (1) that
\begin{eqnarray}
g(S)=\gamma_{M}(G)
\end{eqnarray}
On the other hand,
$S^{'}=a_1a_{n}(\prod\limits_{i=1}^{n-3}a_{i+1}a_{i}^{-1})a_{n-2}^{-1}a_{n}^{-1}$
is the associated surface of $\widetilde{T}^{*}_{\sigma}$, where
$\widetilde{T}^{*}_{\sigma}$ is a joint-tree of $G-v_{2n}$.
Performing a sequence of Transform 4 on $S^{'}$, we have
\begin{eqnarray}
S^{'}  & = & a_1a_{n}(\prod\limits_{i=1}^{n-3}a_{i+1}a_{i}^{-1})a_{n-2}^{-1}a_{n}^{-1} \nonumber \\
 & \sim  & \left\{
        \begin{array}{ll}
        a_1a_{n}a_1^{-1}a_{n}^{-1}(\prod\limits_{i=2}^{n-4}a_{i+1}a_{i}a_{i+1}^{-1}a_{i}^{-1}) & \  \mbox{$n \equiv 0(mod \ 2)$;} \\
        a_1a_{n}a_1^{-1}a_{n}^{-1}(\prod\limits_{i=2}^{n-3}a_{i+1}a_{i}a_{i+1}^{-1}a_{i}^{-1}) & \ \mbox{$n \equiv 1(mod \ 2)$.}  \\
        \end{array}
        \right.
\end{eqnarray}
It can be inferred from (3) that
\begin{eqnarray}
g(S^{'})=\gamma_{M}(G-v_{2n}).
\end{eqnarray}
From (1) and (3) we have
\begin{eqnarray}
g(S)=g(S^{'})+1.
\end{eqnarray}
From (2), (4), and (5) we have $\gamma_{M}(G-v_{2n})$ =
$\gamma_{M}(G)-1$.

According to the above, we can get Theorem 2.1.
$\hspace*{\fill}\Box$

\bigskip

Let $G$ be a connected graph with minimum degree at least 3. The
following algorithm can be used to get the maximum genus of $G$.

\bigskip

{\bf Algorithm I}   \ \  Step 1:  Input $i = 0$, $G_{0}=G$.

Step 2: If there is a 1-$critical$-$vertex$ $v$ in $G_{i}$, then
delete $v$ from $G_{i}$ and go to Step 3. Else, go to Step 4.

Step 3: Deleting all the vertices of degree one and merging all the
vertices of degree two in $G_{i}-v$, we get  a new graph $G_{i+1}$.
Let $i=i+1$, then go back to Step 2.

Step 4: Output $\gamma_{M}(G)= \gamma_{M}(G_{i})+i$.

\bigskip

{\bf Remark }   \ \ Using Algorithm I, the computing of the maximum
genus of $G$ can be reduced to the computing of the maximum genus of
$G_{i}$, which may be much easier than that of $G$.

\bigskip

\noindent {\bf 3. Upper embeddability of  graphs }

\bigskip

An $ear$ of a graph $G$, which is the same as the definition offered
in \cite{wes}, is a path that is maximal with respect to internal
vertices having degree 2 in $G$ and is contained in a cycle in $G$.
An $ear$ $decomposition$ of $G$ is a decomposition $p_0$, \dots,
$p_{k}$ such that $p_0$ is a cycle and $p_{i}$ for $i\geqslant1$ is
an ear of $p_{0}\cup \dots \cup p_{i}$. A $spiral$ $\mathcal
{S}_{m}^{n}$ is the graph which has an ear decomposition $p_0$,
\dots, $p_{n}$ such that $p_0$ is the m-cycle $(v_1v_2\dots v_{m})$,
$p_{i}$ for $1\leqslant i\leqslant m-1$ is the 3-path
$v_{m+2i-2}v_{m+2i-1}v_{m+2i}v_{i}$ which joining $v_{m+2i-2}$ and
$v_{i}$,  and $p_i$ for $i> m-1$ is the 3-path
$v_{m+2i-2}v_{m+2i-1}v_{m+2i}v_{2i-m+1}$ which joining $v_{m+2i-2}$
and $v_{2i-m+1}$. If some edges in $\mathcal {S}_{m}^{n}$ are
replaced by the graph depicted by Fig. 6, then the graph is called
an $extended$-$spiral$, and is denoted by $\mathcal
{\textit{S}}_{m}^{n}$. Obviously, both the vertex  $v_1$ and $v_2$
in Fig. 6 are $\gamma$-vertex. For convenience, a graph of $\mathcal
{S}_{5}^6$ is illustrated by Fig.5, and  Fig. 7 is the graph which
is obtained from $\mathcal {S}_{5}^6$ by replacing the edge
$(v_{13}, v_{14})$ with the graph depicted by Fig. 6.

\begin{footnotesize}

\setlength{\unitlength}{1mm}
\begin{center}
\begin{picture}(100,40)

\put(-8,0){\begin{picture}(10,10)

\put(5,23){\circle{13}}

\put(5,29.3){\circle*{1.5}}

\put(4,26){$v_5$}

\put(3,34){\circle*{1.5}}

\put(-1,35){$v_6$}

\put(-4,28.5){\circle*{1.5}}

\put(-7,30){$v_7$}

\put(-1,30){$p_1$}

\put(-1,25.3){\circle*{1.5}}

\put(0,24){$v_1$}

\qbezier(5,29.3)(6,32.5)(3,34)

\qbezier(-4,28.5)(-3,25)(-1,25.3)

\qbezier(3,34)(-3,34)(-4,28.5)

\qbezier(0.8,18.2)(-11,19)(-4,28.5)

  \put(-6.5,24){\circle*{1.5}}

  \put(-11,24){$v_8$}

  \put(-5,22){$p_2$}

  \put(-4,19){\circle*{1.5}}

  \put(-7,16.5){$v_9$}

\qbezier(-4,19)(2,6)(9.5,18.2)

  \put(-1,14){\circle*{1.5}}

  \put(-5,11.5){$v_{10}$}

  \put(5.5,13.3){\circle*{1.5}}

  \put(4,10.2){$v_{11}$}

  \put(0.8,14.3){$p_3$}

\qbezier(5.5,13.3)(22,14)(11,25.3)

   \put(11.6,14.1){\circle*{1.5}}

   \put(11,11){$v_{12}$}

  \put(11,17){$p_4$}

  \put(16,20){$v_{13}$}

  \put(14.8,20){\circle*{1.5}}

\qbezier(3,34)(18,35)(14.8,20)

   \put(15,27){\circle*{1.5}}

  \put(10,33.2){\circle*{1.5}}

   \put(16,27){$v_{14}$}

  \put(10,29){$p_5$}

  \put(10,34.5){$v_{15}$}

  \put(3,39.5){$v_{16}$}

   \put(-14,34){$v_{17}$}

   \put(2,39){\circle*{1.5}} 

    \put(-10,32){\circle*{1.5}} 

  \qbezier(10,33.2)(6,37)(2,39)  

  \qbezier(-10,32)(-10,30)(-6.5,24) 

  \qbezier(-10,32)(-5,38)(2,39)  

\put(11,25.3){\circle*{1.5}}

   \put(8,23){$v_4$}

  \put(11,17){$p_4$}

  \put(0,19.2){$v_{2}$}

   \put(5.5,18.5){$v_{3}$}

\put(0.8,18.2){\circle*{1.5}}

\put(9.5,18.2){\circle*{1.5}}


\put(-2,2){{\bf Fig. 5.}}

\end{picture}}

\put(45,-5){\begin{picture}(10,10)
\put(-10,24){\circle*{1.5}}

\put(-5,24){\circle*{1.5}}

\put(0,24){\circle*{1.5}}

    \put(4.5,27){\circle*{1.5}}

    \put(4.5,21){\circle*{1.5}}

\put(9,24){\circle*{1.5}}

\put(14,24){\circle*{1.5}}

\put(19,24){\circle*{1.5}}

   \put(4.5,27){\line(-3,-2){5}}

    \put(4.5,27){\line(3,-2){5}}

    \put(4.5,27){\line(0,-1){6}}

   \put(4.5,21){\line(-3,2){5}}

   \put(4.5,21){\line(3,2){5}}

\put(0,24){\line(-1,0){14}}

\put(9,24){\line(1,0){14}}

\qbezier(-10,24)(4,38)(14,24)

\qbezier(-5,24)(5,38)(19,24)

\put(0.8,27.2){{$v_1$}}

\put(5.2,18.5){{$v_2$}}

\put(-6,21){{$v_3$}}

\put(13,21){{$v_4$}}

\put(-1,7){{\bf Fig. 6.}}

\end{picture}}


\put(90,0){\begin{picture}(10,10)


\put(5,23){\circle{13}}

\put(5,29.3){\circle*{1.5}}

\put(4,26){$v_5$}

\put(3,34){\circle*{1.5}}

\put(-1,35){$v_6$}

\put(-4,28.5){\circle*{1.5}}

\put(-7,30){$v_7$}

\put(-1,30){$p_1$}

\put(-1,25.3){\circle*{1.5}}

\put(0,24){$v_1$}

\qbezier(5,29.3)(6,32.5)(3,34)

\qbezier(-4,28.5)(-3,25)(-1,25.3)

\qbezier(3,34)(-3,34)(-4,28.5)

\qbezier(0.8,18.2)(-11,19)(-4,28.5)

  \put(-6.5,24){\circle*{1.5}}

  \put(-11,24){$v_8$}

  \put(-5,22){$p_2$}

  \put(-4,19){\circle*{1.5}}

  \put(-7,16.5){$v_9$}

\qbezier(-4,19)(2,6)(9.5,18.2)

  \put(-1,14){\circle*{1.5}}

  \put(-5,11.5){$v_{10}$}

  \put(5.5,13.3){\circle*{1.5}}

  \put(4,10.2){$v_{11}$}

  \put(0.8,14.3){$p_3$}

\qbezier(5.5,13.3)(22,14)(11,25.3)

   \put(11.6,14.1){\circle*{1.5}}

   \put(11,11){$v_{12}$}

  \put(11,17){$p_4$}

  \put(16,20){$v_{13}$}

  \put(14.8,20){\circle*{1.5}}

\qbezier(3,34)(14,34.5)(15,27)

   \put(15,27){\circle*{1.5}}

  \put(10,33.2){\circle*{1.5}}

   \put(16,27){$v_{14}$}

  \put(10,29){$p_5$}

  \put(10,34.5){$v_{15}$}

  \put(3,39.5){$v_{16}$}

   \put(-14,34){$v_{17}$}

   \put(2,39){\circle*{1.5}} 

    \put(-10,32){\circle*{1.5}} 

  \qbezier(10,33.2)(6,37)(2,39)  

  \qbezier(-10,32)(-10,30)(-6.5,24) 

  \qbezier(-10,32)(-5,38)(2,39)  

\put(11,25.3){\circle*{1.5}}

   \put(8,23){$v_4$}

  \put(11,17){$p_4$}

  \put(0,19.2){$v_{2}$}

   \put(5.5,18.5){$v_{3}$}

\put(0.8,18.2){\circle*{1.5}}

\put(9.5,18.2){\circle*{1.5}}


  \qbezier(24.5,31.3)(35,23)(24.5,15.7)

  \qbezier(24.5,28.3)(35,19)(24.5,12.7)

   \qbezier(15,27)(17,35)(24.5,31.3)

    \qbezier(14.8,20)(19,11)(24.5,12.7)

\put(22,22){\circle*{1.5}}

\put(27,22){\circle*{1.5}}

\put(24.5,31.3){\circle*{1.5}}

\put(24.5,28.3){\circle*{1.5}}

\put(24.5,25.3){\circle*{1.5}}

\put(24.5,25.3){\line(0,1){6}}

\put(24.5,18.7){\line(0,-1){6}}

  \put(24.5,18.7){\circle*{1.5}}

\put(24.5,15.7){\circle*{1.5}}

\put(24.5,12.7){\circle*{1.5}}

 \put(22,22){\line(5,-6){3}}

\put(22,22){\line(5,6){3}}

\put(22,22){\line(1,0){5}}

\put(27,22){\line(-5,6){3}}

\put(27,22){\line(-5,-6){3}}


\put(5,2){{\bf Fig. 7.}}
\end{picture}}

\end{picture}
\end{center}

\end{footnotesize}


{\bf Theorem 3.1} \ \ The graph $\mathcal {S}_{5}^{n}$ is upper
embeddable. Furthermore, $\gamma_{M}(\mathcal {S}_{5}^{n} -
v_{2n+3})$ = $\gamma_{M}(\mathcal {S}_{5}^{n})-1$, $i.e.,$
$v_{2n+3}$ is a 1-$critical$-$vertex$ of $\mathcal {S}_{5}^{n}$.

\bigskip

{\bf Proof} \ \ According to the definition of $\mathcal
{S}_{5}^{n}$, when $n \leq 4$, it is not a hard work to get the
upper embeddability of $\mathcal {S}_{5}^{n}$. So the following 5
cases will be considered.

\medskip

{\bf Case 1:} \ $n=5j$, where $j$ is an integer no less than 1.

\medskip

Without loss of generality, a spanning tree $T$ of  $\mathcal
{S}_{5}^{n}$ can be chosen as $T=T_1\cup T_2$, where $T_1$ is the
path
$v_2v_1v_5v_4v_3\{\prod\limits_{i=1}^{j-1}v_{10i+1}v_{10i}v_{10i-1}v_{10i-2}v_{10i-3}v_{10i-4}
v_{10i+5}v_{10i+4}v_{10i+3}v_{10i+2}\}v_{2n+1}$-
$v_{2n}v_{2n-1}v_{2n-2}v_{2n-3}v_{2n-4}v_{2n+5}v_{2n+4}v_{2n+3}$,
$T_2 = (v_{2n+1}, v_{2n+2})$.  Obviously, the $n+1$ co-tree edges of
$\mathcal {S}_{5}^{n}$ with respect to $T$ are $e_1=(v_2, v_3)$,
$e_2=(v_2, v_9)$, $e_3=(v_1, v_7)$, $\prod\limits_{i=1}^{j-1}
\{e_{5i-1}=(v_{10i-5}, v_{10i-4}), e_{5i}=(v_{10i-6}, v_{10i+3}),
e_{5i+1}=(v_{10i+1}, v_{10i+2}), e_{5i+2}=(v_{10i}, v_{10i+9}),
e_{5i+3}=(v_{10i-2}, v_{10i+7})\}$,  $e_{n-1}=(v_{2n-5}, v_{2n-4})$,
$e_{n}=(v_{2n-6}, v_{2n+3})$, $e_{n+1}=(v_{2n+2}, v_{2n+3})$. Select
such a joint-tree $\widetilde{T}_{\sigma}$ of $\mathcal {S}_{5}^{n}$
which is depicted by Fig.8. After a sequence of Transform 4, the
associated surface $S$ of $\widetilde{T}_{\sigma}$ has the form as
\begin{eqnarray}
S & = & e_1e_2e_1^{-1}e_3e_4e_5
\{\prod\limits_{i=1}^{j-2}e_{5i+1}e_{5i+2}e_{5i-3}^{-1}e_{5i+3}e_{5i-2}^{-1}
e_{5i-1}^{-1}e_{5i+4}e_{5i+5}e_{5i}^{-1}e_{5i+1}^{-1}\} \nonumber\\
& &
e_{n-4}e_{n-3}e_{n-8}^{-1}e_{n-2}e_{n-7}^{-1}e_{n-6}^{-1}e_{n-1}e_{n-5}^{-1}
e_{n-4}^{-1}e_{n-3}^{-1}e_{n-2}^{-1}e_{n-1}^{-1}e_{n+1}e_{n}e_{n+1}^{-1}
e_{n}^{-1} \nonumber \\
& \sim &
\prod\limits_{i=1}^{\lfloor\frac{n+1}{2}\rfloor}e_{i1}e_{i2}e_{i1}^{-1}e_{i2}^{-1},\nonumber
\end{eqnarray}
where $e_{ij}$, $e_{ij}^{-1}$ $\in$ $\{e_1, \dots, e_{n+1},
e_1^{-1}, \dots, e_{n+1}^{-1}\}$;
$i=1,\dots,\lfloor\frac{n+1}{2}\rfloor; j=1, 2.$ Obviously,
$g(S)=\lfloor\frac{n+1}{2}\rfloor$. So, when $n=5j$, $\mathcal
{S}_{5}^{n}$ is upper embeddable.

\medskip


\begin{footnotesize}

\setlength{\unitlength}{1mm}
\begin{center}
\begin{picture}(100,40)

\put(7,-2){\begin{picture}(10,10)

     \put(-20,24){\circle*{1.5}}

\multiput(-20,24)(-1,1){5}{\circle*{0.5}}

\multiput(-20,24)(-1,-1){5}{\circle*{0.5}}

\put(-28,18){{$e_{2}$}}

\put(-28,29){{$e_{1}^{-1}$}}

\put(-21,26.5){{$v_2$}}

  \put(-15,24){\circle*{1.5}}

\put(-16.5,19.5){{$v_1$}}

\multiput(-15,24)(0,1){8}{\circle*{0.5}}

\put(-17,34){{$e_3$}}

  \put(-10,24){\circle*{1.5}}

\put(-11.5,19.5){{$v_5$}}

\multiput(-10,24)(0,1){8}{\circle*{0.5}}

\put(-12,34){{$e_4$}}

   \put(-4.5,24){\circle*{1.5}}

\put(-5,19.5){{$v_4$}}

\multiput(-4.5,24)(0,1){8}{\circle*{0.5}}

\put(-6,34){{$e_5$}}

   \put(1,24){\circle*{1.5}}

\put(-1,26.5){{$v_3$}}

\multiput(1,24)(0,-1){8}{\circle*{0.5}}

\put(0,13){{$e_1$}}

    \put(7,24){\circle*{1.5}}

\put(6,20){{$v_{11}$}}

\multiput(7,24)(0,1){8}{\circle*{0.5}}

\put(6,34){{$e_6$}}

\put(9,27.5){\circle*{0.8}}

\put(12,27.5){\circle*{0.8}}

\put(15,27.5){\circle*{0.8}}

  \put(18,24){\circle*{1.5}}

\multiput(18,24)(0,1){8}{\circle*{0.5}}

\put(15,34){{$e_{n-1}$}}

\put(14,19.5){{$v_{2n-5}$}}

    \put(25,24){\circle*{1.5}}

\multiput(25,24)(0,-1){8}{\circle*{0.5}}

\put(24,13){{$e_{n}^{-1}$}}

\put(20,26.5){{$v_{2n-6}$}}

    \put(32,24){\circle*{1.5}}

\multiput(32,24)(0,1){8}{\circle*{0.5}}

\put(28,34){{$e_{n-5}^{-1}$}}

\put(28,19.5){{$v_{2n-7}$}}

     \put(39,24){\circle*{1.5}}

\put(37,21){{$v_{2n-8}$}}

\multiput(39,24)(0,1){8}{\circle*{0.5}}

\put(38,34){{$e_{n-4}^{-1}$}}

     \put(46,24){\circle*{1.5}}

\put(46,24){\line(3,-5){4}}

\put(49.5,18){\circle*{1.5}}

\put(41,26.5){{$v_{2n+1}$}}

\put(51,16){{$v_{2n+2}$}}

\multiput(49.5,18)(0,-1){7}{\circle*{0.5}}

\put(50,9){{$e_{n+1}^{-1}$}}

  \put(53,24){\circle*{1.5}}

\put(52,26.5){{$v_{2n}$}}

    \put(60,24){\circle*{1.5}}

\multiput(60,24)(0,1){8}{\circle*{0.5}}

\put(56,34){{$e_{n-3}^{-1}$}}

\put(56,19.5){{$v_{2n-1}$}}

  \put(67,24){\circle*{1.5}}

\put(63,26.5){{$v_{2n-2}$}}

   \put(74,24){\circle*{1.5}}

\multiput(74,24)(0,1){8}{\circle*{0.5}}

\put(70,34){{$e_{n-2}^{-1}$}}

\put(69,19){{$v_{2n-3}$}}

    \put(85,24){\circle*{1.5}}

\multiput(85,24)(0,1){8}{\circle*{0.5}}

\put(81,34){{$e_{n-1}^{-1}$}}

\put(80,20){{$v_{2n-4}$}}

  \put(94,24){\circle*{1.5}}

\put(90,27){{$v_{2n+5}$}}

   \put(101,24){\circle*{1.5}}

\put(99,19.5){{$v_{2n+4}$}}

    \put(108,24){\circle*{1.5}}

\multiput(108,24)(0,1){6}{\circle*{0.5}}

\multiput(108,24)(1,-1){4}{\circle*{0.5}}

\put(107,31){{$e_{n+1}$}}

\put(113,19){{$e_{n}$}}

\put(109.5,25){{$v_{2n+3}$}}

\put(-20,24){\line(1,0){128}}

\put(35,2){{\bf Fig. 8.}}
\end{picture}}

\end{picture}
\end{center}
\end{footnotesize}

\medskip

{\bf Case 2:} \ $n=5j+1$, where $j$ is an integer no less than 1.

\medskip

Without loss of generality, select  $T=T_1\cup T_2$ to be a spanning
tree of $\mathcal {S}_{5}^{n}$, where $T_1$ is the path
$v_3v_2v_1\{\prod\limits_{i=1}^{j}v_{10i-3}v_{10i-4}v_{10i-5}v_{10i-6}v_{10i+3}v_{10i+2}
v_{10i+1}v_{10i}v_{10i-1}v_{10i-2}\}v_{2n+5}v_{2n+4}v_{2n+3}$, $T_2
= (v_{2n+1}, v_{2n+2})$.  It is obviously that the $n+1$ co-tree
edges of $\mathcal {S}_{5}^{n}$ with respect to $T$ are $e_1=(v_1,
v_5)$, $e_2=(v_3, v_4)$, $e_3=(v_3, v_{11})$, $e_4=(v_2, v_9)$,
$\prod\limits_{i=1}^{j-1} \{e_{5i}=(v_{10i-3}, v_{10i-2}),
e_{5i+1}=(v_{10i-4}, v_{10i+5}), e_{5i+2}=(v_{10i+3}, v_{10i+4}),
e_{5i+3}=(v_{10i+2}, v_{10i+11}), e_{5i+4}=(v_{10i}, v_{10i+9})\}$,
$e_{n-1}=(v_{2n-5}, v_{2n-4})$, $e_{n}=(v_{2n-6}, v_{2n+3})$,
$e_{n+1}=(v_{2n+2}, v_{2n+3})$. Similar to Case 1, select a joint
tree $\widetilde{T}_{\sigma}$ of $\mathcal {S}_{5}^{n}$. After a
sequence of Transform 4, the associated surface $S$ of
$\widetilde{T}_{\sigma}$ has the form as
\begin{eqnarray}
S & = & e_1e_2e_3e_4e_5e_6e_1^{-1}e_2^{-1}e_7e_8e_3^{-1}e_9
\{\prod\limits_{i=1}^{j-2}e_{5i-1}^{-1}e_{5i}^{-1}e_{5i+5}e_{5i+6}e_{5i+1}^{-1}
e_{5i+2}^{-1}e_{5i+7} \nonumber\\
& &
e_{5i+8}e_{5i+3}^{-1}e_{5i+9}\}e_{n-7}^{-1}e_{n-6}^{-1}e_{n-1}e_{n-5}^{-1}e_{n-4}^{-1}e_{n-3}^{-1}e_{n-2}^{-1}e_{n-1}^{-1}
e_{n+1}e_{n}e_{n+1}^{-1}e_{n}^{-1} \nonumber \\
& \sim &
\prod\limits_{i=1}^{\lfloor\frac{n+1}{2}\rfloor}e_{i1}e_{i2}e_{i1}^{-1}e_{i2}^{-1},
\nonumber
\end{eqnarray}
where $e_{ij}$, $e_{ij}^{-1}$ $\in$ $\{e_1, \dots, e_{n+1},
e_1^{-1}, \dots, e_{n+1}^{-1}\}$;
$i=1,\dots,\lfloor\frac{n+1}{2}\rfloor; j=1, 2.$ Obviously,
$g(S)=\lfloor\frac{n+1}{2}\rfloor$. So, when $n=5j+1$, $\mathcal
{S}_{5}^{n}$ is upper embeddable.

\medskip

{\bf Case 3:} \ $n=5j+2$, where $j$ is an integer no less than 1.

\medskip

Without loss of generality, select a spanning tree of $\mathcal
{S}_{5}^{n}$ to be  $T=T_1\cup T_2$, where $T_1$ is the path
$v_1v_5v_4v_3v_2
\{\prod\limits_{i=1}^{j}v_{10i-1}v_{10i-2}v_{10i-3}v_{10i-4}v_{10i+5}v_{10i+4}
v_{10i+3}v_{10i+2}v_{10i+1}v_{10i}\}
 v_{2n+5}v_{2n+4}v_{2n+3}$, $T_2
= (v_{2n+1}, v_{2n+2})$.  It is obviously that the $n+1$ co-tree
edges of $\mathcal {S}_{5}^{n}$ with respect to $T$ are $e_1=(v_1,
v_2)$, $e_2=(v_1, v_7)$, $e_3=(v_5, v_6)$, $e_4=(v_4, v_{13})$,
$e_5=(v_3, v_{11})$, $\prod\limits_{i=1}^{j-1}
\{e_{5i+1}=(v_{10i-1}, v_{10i}), e_{5i+2}=(v_{10i-2}, v_{10i+7}),
e_{5i+3}=(v_{10i+5}, v_{10i+6}), e_{5i+4}=(v_{10i+4}, v_{10i+13}),
e_{5i+5}=(v_{10i+2}, v_{10i+11})\}$, $e_{n-1}=(v_{2n-5}, v_{2n-4})$,
$e_{n}=(v_{2n-6}, v_{2n+3})$, $e_{n+1}=(v_{2n+2}, v_{2n+3})$.
Similar to Case 1, select a joint-tree $\widetilde{T}_{\sigma}$ of
$\mathcal {S}_{5}^{n}$. After a sequence of Transform 4, the
associated surface $S$ of $\widetilde{T}_{\sigma}$ has the form as
\begin{eqnarray}
S & = &
e_1e_2e_3e_4e_5e_1^{-1}e_6e_7e_2^{-1}e_3^{-1}e_8e_9e_4^{-1}e_{10}
\{\prod\limits_{i=1}^{j-2}e_{5i}^{-1}e_{5i+1}^{-1}e_{5i+6}e_{5i+7}e_{5i+2}^{-1}
e_{5i+3}^{-1}e_{5i+8} \nonumber\\
& &
e_{5i+9}e_{5i+4}^{-1}e_{5i+10}\}e_{n-7}^{-1}e_{n-6}^{-1}e_{n-1}
e_{n-5}^{-1}e_{n-4}^{-1}e_{n-3}^{-1}e_{n-2}^{-1}e_{n-1}^{-1}
e_{n+1}e_{n}e_{n+1}^{-1}e_{n}^{-1} \nonumber \\
& \sim &
\prod\limits_{i=1}^{\lfloor\frac{n+1}{2}\rfloor}e_{i1}e_{i2}e_{i1}^{-1}e_{i2}^{-1},
\nonumber
\end{eqnarray}
where $e_{ij}$, $e_{ij}^{-1}$ $\in$ $\{e_1, \dots, e_{n+1},
e_1^{-1}, \dots, e_{n+1}^{-1}\}$;
$i=1,\dots,\lfloor\frac{n+1}{2}\rfloor; j=1, 2.$ Obviously,
$g(S)=\lfloor\frac{n+1}{2}\rfloor$. So, when $n=5j+2$, $\mathcal
{S}_{5}^{n}$ is upper embeddable.

\medskip

{\bf Case 4:} \ $n=5j+3$, where $j$ is an integer no less than 1.

\medskip

Without loss of generality, a spanning tree $T$ of  $\mathcal
{S}_{5}^{n}$ can be chosen as $T=T_1\cup T_2$, where $T_1$ is the
path $v_2v_1v_7v_6v_5v_4v_3
\{\prod\limits_{i=1}^{j}v_{10i+1}v_{10i}v_{10i-1}v_{10i-2}v_{10i+7}v_{10i+6}
v_{10i+5}v_{10i+4}v_{10i+3}$-$v_{10i+2}\}
 v_{2n+5}v_{2n+4}v_{2n+3}$, $T_2
= (v_{2n+1}, v_{2n+2})$.  It is obviously that the $n+1$ co-tree
edges of $\mathcal {S}_{5}^{n}$ with respect to $T$ are $e_1=(v_1,
v_5)$, $e_2=(v_2, v_3)$, $e_3=(v_2, v_9)$, $\prod\limits_{i=1}^{j}
\{e_{5i-1}=(v_{10i-3}, v_{10i-2}), e_{5i}=(v_{10i-4}, v_{10i+5}),
e_{5i+1}=(v_{10i-6}, v_{10i+3}), e_{5i+2}=(v_{10i+1}, v_{10i+2}),
e_{5i+3}=(v_{10i}, v_{10i+9})\}$, $e_{n+1}=(v_{2n+2}, v_{2n+3})$.
Similar to Case 1, select a joint-tree $\widetilde{T}_{\sigma}$ of
$\mathcal {S}_{5}^{n}$. After a sequence of Transform 4, the
associated surface $S$ of $\widetilde{T}_{\sigma}$ has the form as
\begin{eqnarray}
S & = & e_1e_2e_3e_4e_5e_1^{-1}e_6e_2^{-1}
\{\prod\limits_{i=1}^{j-1}e_{5i+2}e_{5i+3}e_{5i-2}^{-1}e_{5i-1}^{-1}e_{5i+4}e_{5i+5}
e_{5i}^{-1}e_{5i+6} \nonumber\\
& & e_{5i+1}^{-1}e_{5i+2}^{-1}\}e_{n-1}e_{n-5}^{-1}
e_{n-4}^{-1}e_{n-3}^{-1}e_{n-2}^{-1}e_{n-1}^{-1}
e_{n+1}e_{n}e_{n+1}^{-1}e_{n}^{-1} \nonumber \\
& \sim &
\prod\limits_{i=1}^{\lfloor\frac{n+1}{2}\rfloor}e_{i1}e_{i2}e_{i1}^{-1}e_{i2}^{-1},
\nonumber
\end{eqnarray}
where $e_{ij}$, $e_{ij}^{-1}$ $\in$ $\{e_1, \dots, e_{n+1},
e_1^{-1}, \dots, e_{n+1}^{-1}\}$;
$i=1,\dots,\lfloor\frac{n+1}{2}\rfloor; j=1, 2.$ Obviously,
$g(S)=\lfloor\frac{n+1}{2}\rfloor$. So, when $n=5j+3$, $\mathcal
{S}_{5}^{n}$ is upper embeddable.

\medskip

{\bf Case 5:} \ $n=5j+4$, where $j$ is an integer no less than 1.

\medskip

Without loss of generality, a spanning tree $T$ of  $\mathcal
{S}_{5}^{n}$ can be chosen as $T=T_1\cup T_2\cup T_3$, where $T_1$
is the path $v_1v_2
\{\prod\limits_{i=1}^{j}v_{10i-1}v_{10i-2}v_{10i-3}v_{10i-4}v_{10i-5}v_{10i-6}
v_{10i+3}v_{10i+2}v_{10i+1}v_{10i}\}
 v_{2n+1}$-$v_{2n}v_{2n-1}v_{2n-2}v_{2n-3}v_{2n-4}v_{2n+5}v_{2n+4}v_{2n+3}$,
 $T_2
= (v_2, v_3)$,
 $T_3
= (v_{2n+1}, v_{2n+2})$.  It is obviously that the $n+1$ co-tree
edges of $\mathcal {S}_{5}^{n}$ with respect to $T$ are $e_1=(v_1,
v_5)$, $e_2=(v_1, v_7)$, $e_3=(v_3, v_4)$,  $e_4=(v_3, v_{11})$,
$\prod\limits_{i=1}^{j} \{e_{5i}=(v_{10i-1}, v_{10i}),
e_{5i+1}=(v_{10i-2}, v_{10i+7}), e_{5i+2}=(v_{10i-4}, v_{10i+5}),
e_{5i+3}=(v_{10i+3}, v_{10i+4}), e_{5i+4}=(v_{10i+2},
v_{10i+11})\}$, $e_{n+1}=(v_{2n+2}, v_{2n+3})$. Similar to Case 1,
select a joint-tree $\widetilde{T}_{\sigma}$ of $\mathcal
{S}_{5}^{n}$. After a sequence of Transform 4, the associated
surface $S$ of $\widetilde{T}_{\sigma}$ has the form as
\begin{eqnarray}
S & = & e_2e_1e_3e_4e_5e_6e_2^{-1}e_7e_1^{-1}e_3^{-1}
\{\prod\limits_{i=1}^{j-1}e_{5i+3}e_{5i+4}e_{5i-1}^{-1}e_{5i}^{-1}e_{5i+5}e_{5i+6}
e_{5i+1}^{-1} \nonumber\\
& & e_{5i+7}e_{5i+2}^{-1}e_{5i+3}^{-1}\}e_{n-1}e_{n-5}^{-1}
e_{n-4}^{-1}e_{n-3}^{-1}e_{n-2}^{-1}e_{n-1}^{-1}
e_{n+1}e_{n}e_{n+1}^{-1}e_{n}^{-1} \nonumber \\
& \sim &
\prod\limits_{i=1}^{\lfloor\frac{n+1}{2}\rfloor}e_{i1}e_{i2}e_{i1}^{-1}e_{i2}^{-1},
\nonumber
\end{eqnarray}
where $e_{ij}$, $e_{ij}^{-1}$ $\in$ $\{e_1, \dots, e_{n+1},
e_1^{-1}, \dots, e_{n+1}^{-1}\}$;
$i=1,\dots,\lfloor\frac{n+1}{2}\rfloor; j=1, 2.$ Obviously,
$g(S)=\lfloor\frac{n+1}{2}\rfloor$. So, when $n=5j+4$, $\mathcal
{S}_{5}^{n}$ is upper embeddable.

From the Case 1-5, the upper embeddability of  $\mathcal
{S}_{5}^{n}$ can be obtained.

Similar to the Case 1-5, for each $n\geq 5$, there exists a
joint-tree $\widetilde{T}^{*}_{\sigma}$ of $\mathcal
{S}_{5}^{n}-v_{2n+3}$ such that its associated surface is
$S^{'}=S-\{e_{n+1}e_{n}e_{n+1}^{-1}e_{n}^{-1}\}$. It is obvious that
$S^{'}$ is the surface into which the embedding of $\mathcal
{S}_{5}^{n}-v_{2n+3}$ is the maximum genus embedding. Furthermore,
$g(S^{'})=g(S)-1$, $i.e.,$ $\gamma_{M}(\mathcal {S}_{5}^{n} -
v_{2n+3})$ = $\gamma_{M}(\mathcal {S}_{5}^{n})-1$. So, $v_{2n+3}$ is
a 1-$critical$-$vertex$ of $\mathcal {S}_{5}^{n}$. $\hspace*{\fill}
\Box$

\bigskip

Similar to the proof of Theorem 3.1, we can get the following
theorem.

\medskip

{\bf Theorem 3.2} \ \ The graph $\mathcal {S}_{m}^{n}$ is upper
embeddable. Furthermore, $\gamma_{M}(\mathcal {S}_{m}^{n} -
v_{m+2n-2})$ = $\gamma_{M}(\mathcal {S}_{m}^{n})-1$, $i.e.,$
$v_{m+2n-2}$ is a 1-$critical$-$vertex$ of $\mathcal {S}_{m}^{n}$.

\bigskip

{\bf Corollary 3.1} \ \ Let $G$ be a graph with minimum degree at
least three. If $G$, through a sequence of vertex-splitting
operations, can be turned into a  $spiral$ $\mathcal {S}_{m}^{n}$,
then $G$ is upper embeddable.

\bigskip

{\bf Proof} \ \ According to Lemma 1.3, Theorem 3.2, and the upper
embeddability of graphs, Corollary 3.1 can be obtained.
$\hspace*{\fill} \Box$

\bigskip

In the following, we will offer an  algorithm to obtain the maximum
genus of the $extended$-$spiral$ $\mathcal {\textit{S}}_{m}^{n}$.

\bigskip

{\bf Algorithm II} \ \ Step 1: Input $i=0$ and  $j=0$.  Let $G_0$ be
the $extended$-$spiral$ $\mathcal {\textit{S}}_{m}^{n}$.

Step 2: If there is a $\gamma$-vertex $v$ in $G_{i}$, then delete
$v$ from $G_{i}$, and go to Step 3. Else, go to Step 4.

Step 3: Deleting all the vertices of degree one and merging some
vertices of degree two in $G_{i}-v$, we get  a new graph $G_{i+1}$.
Let $i=i+1$. If $G_{i}$ is a $spiral$ $\mathcal {S}_{m}^{n}$, then
go to Step 4. Else, go back to Step 2.

Step 4: Let $G_{i+j}$ be the $spiral$ $\mathcal {S}_{m}^{n}$.
Deleting $v_{m+2n-2}$ from $\mathcal {S}_{m}^{n}$, we will get a new
graph $G_{i+j+1}$, (obviously, $G_{i+j+1}$ is either a $spiral$
$\mathcal {S}_{m}^{n-2}$ or a $cactus$).

Step 5: If $G_{i+j+1}$ is a $cactus$, then go to Step 6. Else, Let
$n=n-2$, $j=j+1$ and go back to Step 4.

Step 6: Output $\gamma_{M}(\mathcal {\textit{S}}_{m}^{n})= i+j+1$.

\bigskip

{\bf Remark } \ 1. \  In the graph $G$ depicted by Fig. 6, after
deleting a $\gamma$-vertex $v_1$ (or $v_2$) from $G$, the vertex
$v_3$ (or $v_4$) is still a $\gamma$-vertex of the remaining graph.

2.  \   From Algorithm II we can get that the $extended$-$spiral$
$\mathcal {\textit{S}}_{m}^{n}$ is upper embeddable.

 
\end{document}